\documentclass[11pt ]{article}

\usepackage{latexsym,amssymb}
\usepackage[applemac]{inputenc}

\newif\ifpdf
\ifx\pdfoutput\undefined
\pdffalse % on n'utilise pas pdflatex
\else
\pdfoutput=1 %on utilise  pdflatex
\pdftrue
\fi
\ifpdf
\usepackage[pdftex]{graphicx}
\else
\usepackage[dvips]{graphicx}
\fi
 \newcommand{\Section}{\setcounter{equation}{0} \section}

\newcommand{\N} {\mathbb{N}}             %
\newcommand{\R} {\mathbb{R}}     
\newcommand{\eps} {\varepsilon}          %
  \newcommand{\preuve}{\noindent\textit{ Proof -~}}
 \newcommand{\findemo}{\hfill $\Box$}
\newcommand{\cP} {\mathcal{P}}
\newcommand{\cS} {\mathcal{S}}
\newcommand{\wcS}{\widetilde{\mathcal{S}}}
\newcommand{\wcP}{\widetilde{\mathcal{P}}}
\newtheorem{theo}{Theorem}[section]
\newtheorem{prop}{Proposition}[section]
\newtheorem{cor}{Corollary}[section]

\newtheorem{lemme}{Lemma}[section]

\newtheorem{rem}{Remark}[section]
  
\begin{document}
\ifpdf
\DeclareGraphicsExtensions{.pdf,.jpg}
\else
\DeclareGraphicsExtensions{.eps,.jpg}
\fi
\thispagestyle{empty}
   \begin{center}
{\Large{\textbf{ A regularization method for  ill-posed bilevel optimization problems}}}

M. Bergounioux \& M. Haddou \\
MAPMO-UMR 6628\\
Universit\'e d'orl\'eans - BP 6759 \\
45067 Orl\'eans cedex 2 \\
maitine.bergounioux, mounir.haddou@univ-orleans.fr

\today
\end{center}

\begin{abstract}
We present a regularization method to approach a solution of the pessimistic formulation of ill -posed bilevel problems . This allows to overcome the difficulty arising from the non uniqueness of the lower level problems solutions and responses. We prove existence of approximated solutions, give convergence result using Hoffman-like assumptions. We end with objective  value 
error estimates. 
\end{abstract}

 \Section{Introduction}
 Bilevel programming problems are of growing interest both from theoretical 
and   practical points of view. These models  are used  in various  applications, such as economic planning, network design, and so on...
This large class of important and strategic economic problems   can be viewed 
 as  static  noncooperative asymmetric games. Two players seek to optimize their individual objective functions.
 The first player (the leader)  must take  into account the reaction
 (or any possible reaction when non unique) of  the second player (the follower).
 Such a problem  can be ill posed since  the lower level ( follower's  problem ) may  have many solutions and  responses for every  (or some) fixed leader's variables.
In the so-called ``optimistic case '', many  optimal reactions of the follower are possible and the follower is assumed to choose in favor of the leader. In this case,  the upper level  problem can be modelled using a bilevel formulation.
These programs are quite difficult nonconvex optimization problems. Several theoretical results and heuristics or approximation techniques can be found in the recent literature \cite{DEMPE, LM1,LM2,BLMS,LMS,MS,NBL}. In some of these works, strong assumptions are made to simplify the model. The solution  of the lower level is supposed to be unique or   if they are many, they provide the same (unique) upper level objective value.
 
 The lower level is replaced by the equivalent first order optimality conditions and can be viewed as an equilibrium problem. Most of time, the complementarity part of these optimality conditions is smoothed or penalized using different techniques.\\
 In our approach we consider the realistic situation where different  reactions of the follower are possible. 
There are multiple responses and  we consider the so-called ``pessimistic''  formulation of the asymmetric game.  
This can be interpreted as a kind of non-cooperative asymmetric game

Throughout this paper, we shall consider the general bi-level problem:
 $$(\cP)~\left \{ \begin{array}{l}
 \max f(y,x) \\
 y \in K~,~x \in \cS(y)~, \end{array}\right .$$
 where $K$ and $C$ are non empty  convex, closed, bounded subsets of $\R^n$  and
 \begin{equation} \label{Sy}
 \cS(y) = \mbox{argmin} \,\{\,h(y,z)~| z \in C~\},
 \end{equation}
  $f$ and $h $ are smooth functions from $\R^n\times \R^m$ to $\R$. 
  Moreover, for every $y\in K$, $f(y,\cdot )$ and $h(y,\cdot )$ are convex  and $f$ takes only positive values ( assumptions will be made more precise later). 
 \begin{rem} 
  Since the upper-level objective function $f$ has to be maximized, one can suppose it (without loss of generality) to be positive. Indeed, $f$ can be replaced by  \\
  $\left (\max \{  f-f(y_0,x_0), 0\}\right )^2$  for some $y_0\in K$ and $x_0 \in \cS(y_0)$ . 
    \end{rem}
  As mentionned before, the main difficulty comes from the fact that the  cost ``function''
 $ f(y,x)~, ~x \in \cS(y)$ 
  can be a \textbf{multivalued application}  whenever the lower level set is not unique and distinct
  solutions yield distinct upper-level objective function values.  
  In addition, it is not clear that 
  $ f(y,x) = f(y,\tilde x) $ for any $x,\tilde x \in \cS(y)$. Therefore, it is difficult to compute the solutions
   (if there are any).

The remaining of the paper is organized as follows. We introduce  the regularization method and give an existence result in next section. Section 3 is devoted to an asymptotic analysis: we prove that the cluster points of  solutions to the penalized problems are solutions of a three-level limit problem corresponding to the `` pessimistic''  formulation of the considered  ill-posed bilevel problems . We  give some error estimates  in the last section.

 %%%%%%%%%%%%%%%%%%%%%%
\Section{The penalized problem}
%%%%%%%%%%%%%%%%%%%%%%%  
   We would like to let the upper level objective function  single valued. So we are going to use a penalization process that allow to compute approximate  solutions  more easily. More precisely, $\eps >0$ being given, we consider the following penalized problem
  $$( \cP_\eps) ~\left \{ \begin{array}{l}
 \max f(y,x) \\
 y \in K~,~x \in \cS_\eps(y)~, \end{array}\right .$$
where 
\begin{equation} \label{Seps}
 \cS_\eps (y) = \mbox{argmin}\, \{\,  h_\eps (y,z)~| z \in C~\},
 \end{equation}
where
\begin{equation} \label{geps}
h_\eps (y) = h(y,z) + \eps f^2(y,z)~. 
 \end{equation}
For each nonnegative $\eps$, the bi-level problem $( \cP_\eps) $ is well posed . Furthermore, under  some general and non restrictive assumptions on $f$ and $h$ we will prove that the upper level function is single valued and continuous with respect to the leader variables $y$. 

This regularization technique makes some selection property on the solutions of the lower level problem
which is easy to characterize and have an explicit and simple  economic interpretation. 
In almost all other regularization methods, the lower level is replaced by its optimality conditions. The bi-level problem is then considered as a mathematical program with equilibrium constraints.  The ``hard''  part of these constraints ( namely the complementarity conditions ) is  then smoothed or penalized.  
In fact  these methods make also some selection ( the generated sequences converge to the analytic center when using smoothing methods  or the least norm center ) on the solution set of the lower level but these selections do not have any economic interpretation since they have no link to the objective function of the upper level.  Moreover, convergence results need more restrictive assumptions.  \\
%From now  on, we assume that 
%\begin{equation}
%f \mbox{ and } h \mbox{ are  continous with respect to both variables } x \mbox{ and } y~.
%\end{equation}
For convenience of the reader, we first give  or recall some direct and classical results. These results will be useful for forthcoming developpements.
\begin{lemme}\label{lem1} For any $\eps >0$, the lower-level problem 
$$ \mathcal{Q}_{\eps, y} = \left \{ \begin{array}{l}
\min h_\eps (y,z) \\
z \in C~,
\end{array} \right . $$
admits (at least) a solution so that $ \cS_\eps (y)\neq \emptyset$. Moreover, there exists a constant $ \kappa_y \in \R  $ such that 
$$ \forall x \in  \cS_\eps (y) \qquad   f(y,x)  = \kappa_y ~.$$
\end{lemme}
\preuve The existence of a solution to $ \mathcal{Q}_{\eps, y} $ is obvious since $C$ is bounded and $f,~h$ are continuous.  Moreover $ \mathcal{Q}_{\eps, y} $ may be written as follows
$$ \mathcal{Q}^*_{\eps, y} ~~ \left \{ \begin{array}{l}
\min  h (y,z) + \eps t^2\\
f(y,z)-t= 0~,\\
z \in C~,
\end{array} \right . $$
Since $f$ takes only positive values, the equality constraint can be obviously replaced by an inequality in this minimization problem
 $$ \mathcal{Q}^*_{\eps, y} ~~ \left \{ \begin{array}{l}
\min  h (y,z) + \eps t^2\\
f(y,z)-t\le 0~,\\
z \in C~.
\end{array} \right . $$
$\mathcal{Q}^*_{\eps, y} $ is a convex problem and the cost function  is strictly convex with respect to $t$.\\
This simple observation proves that  the optimal value of $t$ is unique and completes the proof.
 \findemo
 \begin{lemme}\label{lem2} Let be $\eps >0$ fixed. The multi-application $\cS_\eps$ is lower semi-continuous in the following sense :
 if $y_k \to y$ and $x_k \in \cS_\eps(y_k)$ then $ x_k \to x \in  \cS_\eps(y)$ (up to a subsequence).
\end{lemme}
\preuve Let be $x_k \in \cS_\eps(y_k)\subset C$. As $C$ is bounded, then $(x_k )$ is bounded as well and converges to some $x$ (up to a subsequence). As $x_k \in \cS_\eps(y_k)$ we get
$$ \forall z \in C \qquad h(y_k,x_k) +\eps  f^2(y_k,x_k) \le 
 h(y_k,z) +\eps f^2(y_k,z) ~.$$
 As $f$ and $h$ are continuous with respect to $y$ and $x$ we obtain
 $$ \forall z \in C \qquad h(y,x ) +\eps  f^2(y,x) \le 
 h(y,z) +\eps  f^2(y,z)~,$$
 that is $x \in \cS_\eps(y)$.
\findemo 
\begin{lemme}\label{lem3} Let be $\eps >0$ fixed. The cost function 
$$ v_\eps : y \mapsto \{ f(y,x)~|~x \in  \cS_\eps(y)\,\}$$
is single-valued and continuous. \end{lemme}
\preuve We see that the function $v_\eps$ is single valued, with Lemma \ref{lem1}. Let us prove the continuity: let be $(y_k)$ a sequence that converges to some $y$. Then 
$v_\eps (y_k) = f(y_k, x_k) $ where $x_k \in \cS_\eps(y_k)$. Lemma \ref{lem2} yields that $x_k$ converges (up to a subsequence) to $x \in \cS_\eps(y)$. As $f$ is continuous with respect to $y$ and $x$ we get 
$$ v_\eps (y_k) =  f(y_k, x_k)  \to  f(y, x) = v_\eps(y)~.$$
~
\findemo
\\
We may now give an existence result :
\begin{theo} For any $\eps >0$, problem $( \cP_\eps) $ admits at least an optimal solution 
$y_\eps$.
\end{theo}
\preuve As $v_\eps$ is continuous and $K$ is bounded, the result follows.\findemo
%%%%%%%%%%%%%%%%%%%%%%~
\Section{Asymptotic results }
%%%%%%%%%%%%%%%%%%%%%%%
\subsection{ A convergence result for the solutions of  $(\cP_\eps)$}
 In this subsection, we study the behaviour of solutions of $( \cP_\eps) $ as $\eps$ goes to 0.  
 First, we introduce some notations:
\begin{equation}\label{stildey}\wcS(y) = \mbox{argmin} \{   f^2(y,z)~| z \in \cS(y)~\},
\end{equation}
where 
$ \cS (y) $ is given by (\ref{Sy})
and 
\begin{equation}\label{ptilde}(\wcP)~\left \{ \begin{array}{l}
 \max f(y,x) \\
 y \in K~,~x \in \wcS(y)~, \end{array}\right .\end{equation}
  Note that problem $(\wcP)$ is a three-level problem that can be written in an extended way as follows:
  $$(\wcP)~\left \{ \begin{array}{l}
 \max f(y,x) \\
 y \in K\\
 x \in  \mbox{argmin} \left \{   f^2(y,z)~| z \in \mbox{argmin} \,\{\,h(y,w)~| w \in C~\}\right \},\end{array}\right .$$

\begin{lemme}\label{lem11}  $x \in \wcS(y) $  is equivalent to  
 $$x\in \cS(y) \mbox{ and }
 \forall z \in C \mbox{ such that } h(y,z)=h(y,x),\qquad f^2(y,z) \ge f^2(y,x)~.$$
 \end{lemme}
 \preuve Assume that $z$ satisfies $ h(y,z)=h(y,x) $ with $x \in $ argmin  $\{\, h(y,t) ~|~t\in C\, \}$. Then $ z \in $ argmin  $\{\, h(y,t) ~|~t\in C\, \}$. 
 \findemo
 \begin{lemme}\label{lem12}  Let $y$ be fixed. If $x_\eps \in \cS_\eps$ converges to some $\bar x$, then $\bar x \in  \wcS(y) $\end{lemme}
 \preuve Assume $x_\eps \in \cS_\eps$ and $x_\eps \to \bar x$ as $\eps \to 0$. 
For every $z\in C$ we get 
 $$  h(y,x_\eps) + \eps f^2(y, x_\eps) \le   h(y,z) + \eps f^2(y,z)~.$$
When $\eps  \to 0$,  as the functions are continuous we obtain 
 $$\forall z \in C \qquad   h(y,\bar x ) \le   h(y,z) ~,$$
that is $\bar x \in \cS(y)$.
\\
Let be $\tilde x \in C $ such that $h(y,\tilde x) = h(y, \bar x)$. Then 
$$ \begin{array}{rcll} h(y,x_\eps) + \eps f^2(y, x_\eps) &\le&  h(y,\tilde x) + \eps f^2(y,\tilde x)&\mbox{since } \tilde x \in C\\
&\le & h(y,\bar x) + \eps f^2(y,\tilde x) &\mbox{since }h(y,\tilde x) = h(y, \bar x)\\
&\le & h(y,x_\eps ) + \eps f^2(y,\tilde x) &\mbox{since } x_\eps \in C \mbox{ and } \bar x \in \cS(y)~.
\end{array}$$
Therefore 
$$ \forall \tilde x \in C  \mbox{  such that } h(y,\tilde x) = h(y, \bar x), ~f^2(y, x_\eps)\le f^2(y,\tilde x) ~.$$
Passing to the limit with the continuity of $f$ gives 
$$ \forall \tilde x \in C  \mbox{  such that } h(y,\tilde x) = h(y, \bar x), ~f^2(y, \bar x)\le f^2(y,\tilde x) ~.$$
With Lemma \ref{lem11} we conclude that  $ x \in \wcS(y) $.   \findemo

To establish the main convergence result of this work, we will use some technical but not so  very restrictive assumption.\\
Let us  set
\begin{equation}\label{alpheps}
\alpha_\eps = h(y_\eps, x_\eps) + o(\eps)~,
\end{equation} 
and 
\begin{equation}\label{Leps}
\Lambda_\eps =\{ \, x \in C \, |\,  h(y_\eps, x )  \le \alpha_\eps\, \}~.
\end{equation} 
Assume we can find $\sigma_o >0$ and $\eps_o >0$ such that 
\begin{equation}\label{hyp1} \forall \eps  \le \eps_o~~ \inf_{h(y_\eps, x) = \alpha_\eps } |\nabla_x h (y_\eps, x)|\ge \sigma_o~,\end{equation}
This assumption does not seem quite natural at a first glimpse. In fact it a Hoffman- inequality type assumption which is more or less standard in this context.  The proof of next theorem, and especially the  proof of Lemma \ref{lemstar} below will make this hypothesis clear. 
\begin{theo}\label{asym1} Assume condition (\ref{hyp1}) is verified  and let  $ y_\eps $ an optimal solution to $(\cP_\eps)$. Then $y_\eps$ converges to some $\bar y$ (up to a subsequence) and $\bar y $ is an optimal solution to $(\wcP)$. 
\end{theo}
\preuve Let  $ y_\eps $ an optimal solution to $(\cP_\eps)$. Then $y_\eps \in K$ which is bounded. So (extracting a subsequence) we may assert that $y_\eps$ converges to $\bar y$. As $K$ is closed then $\bar y \in K$.  As $ y_\eps $ is  an optimal solution to $(\cP_\eps)$ we have
\begin{equation}\label{eq31}
\forall \tilde y  \in K~, \forall \tilde x_\eps \in \cS_\eps(\tilde y) \qquad  f(y_\eps, x _\eps ) \ge f(\tilde y, \tilde x_\eps)
\end{equation}
where $x_\eps \in \cS_\eps(y_\eps) $. Note that $ \tilde x_\eps \in \cS_\eps(\tilde y) $ implies that $\tilde x_\eps \in C$. So $\tilde x_\eps$ is bounded and converges to $\tilde x $ (up to a subsequence) with $\tilde x \in \wcS(\tilde y)$ (Lemma \ref{lem12}). 
\\ Passing to the limit in (\ref{eq31}) gives 
$$ \forall \tilde y  \in K ~, \exists \tilde x \in \wcS(\tilde y) \mbox{ such that }  f(\bar y, \bar x) \ge f(\tilde y, \tilde x )~,$$
where $\bar x $ is the limit (of a subsequence) of $x_\eps$.  Now we   need  the  following result to achieve the proof :
%**********************
\begin{lemme}\label{lemstar} Assume that   (\ref{hyp1})  is satisfied and let $(y_\eps, x_\eps\in \cS_\eps(y_\eps)) $ converging to $(\bar y, \bar x)$. Then $\bar x \in \wcS(\bar y)$.
\end{lemme} 
%**********************
Thanks to the definition of $ \wcS(\tilde y)$  we note that $f(\tilde y, \cdot)$ is constant on $ \wcS(\tilde y)$, namely
$$\forall z \in   \wcS(\tilde y)  \qquad f(\tilde y, z)= f(\tilde y, \tilde x )~.$$
Finally
$$ \forall \tilde y  \in K,~\forall \tilde x \in \wcS(\tilde y) \qquad  f(\bar y, \bar x) \ge f(\tilde y, \tilde x )~,$$
with $\bar x \in  \wcS(\bar y)$. This means that $\bar y$ is an optimal solution to $(\wcP)$. 
\\
$\bullet$ It remains to prove  Lemma \ref{lemstar}.  

 Let $y_\eps $ converging to $\bar y$ and  $x_\eps\in \cS_\eps(y_\eps)$. As $x_\eps\in  C$ (bounded) one may extract a subsequence converging 
to $\bar x$.   We are going to prove that $\bar x \in \wcS(\bar y)$.\\

We first prove that $\bar x \in \cS(\bar y)$.  As $x_\eps \in \cS_\eps(y_\eps)$ we have 
\begin{equation}\label{eq32}
\forall z \in C  \qquad  h(y_\eps, x_\eps  ) + \eps \, f^2( y_\eps, x_\eps) \le  h(y_\eps, z ) + \eps \, f^2( y_\eps,z) ~;\end{equation} 
as $f$ and $g$ are continuous, passing to the limit gives 
$$\forall z \in C  \qquad  h(\bar y , \bar x  ) \le  h(\bar y , z ) ~,$$ that is 
$\bar x \in \cS(\bar y)$. \\
Let $\tilde x \in \cS(\bar y)$.  Suppose for a while that $\exists \tilde \eps $ such that 
\begin{equation}\label{eq33}\forall \eps \le \tilde \eps \qquad \tilde x \in \Lambda_\eps~.\end{equation}
We get 
$$ h( y_\eps, \tilde x ) \leq h(y_\eps, x_\eps) + o(\eps)~;$$
with relation (\ref{eq32})  this gives 
\begin{equation}\label{eq33b}\forall z \in C  \qquad  h(y_\eps, \tilde x ) + \eps \, f^2( y_\eps, x_\eps) \le  h(y_\eps, z ) + \eps \, f^2( y_\eps,z)  + o(\eps)~.
\end{equation}
 As $\tilde x \in C$ relation (\ref{eq32}) yields as well
$$  h(y_\eps, x_\eps  ) + \eps \, f^2( y_\eps, x_\eps) \le  h(y_\eps, \tilde x ) + \eps \, f^2( y_\eps,\tilde x) ~.$$ 
Adding these two relations gives
 \begin{equation}\label{eq34} \forall z \in C  \qquad h(y_\eps, x_\eps  ) + 2 \eps \, f^2( y_\eps, x_\eps) \le  h(y_\eps, z ) + \eps \, f^2( y_\eps,z) + \eps \, f^2( y_\eps,\tilde x)   + o(\eps)~; \end{equation}
the choice of $z= x_\eps$ implies 
$$   \eps \, f^2( y_\eps, x_\eps) \le    \eps \, f^2( y_\eps,\tilde x)   + o(\eps)~,$$ 
that is 
$$f^2( y_\eps, x_\eps) \le     f^2( y_\eps,\tilde x)   + \frac{o(\eps)}{\eps}~.$$ 
Passing to the limit gives finally
$$ \forall\tilde x \in \cS(\bar y)\qquad f^2(\bar y , \bar x ) \le     f^2( \bar y ,\tilde x)   ~.$$ 
This means that $\bar x \in \wcS(\bar y)$. 

Unfortunately,  there is no reason for ``assumption '' (\ref{eq33})  to be  satisfied and we must get rid of it. We are going  to adapt the previous proof (we gave the main ideas). If 
$\tilde x \notin \Lambda_\eps$ then we perform a projection: we call $\tilde x_\eps $ the projection of $\tilde x$ on $\Lambda_\eps$.  We are going to show that $\tilde x_\eps $ converges to $\tilde x$. \\
As $\tilde x \notin \Lambda_\eps$  we get  $ \alpha_\eps < h(y_\eps, \tilde x)$ . Let us call $\sigma_{\alpha_\eps} (h) $ the following real number 
\begin{equation}\label{hofman}
\sigma_{\alpha_\eps} (h)  = \inf_{x\in [  \alpha_\eps < h(y_\eps, \cdot )] } \frac{h(y_\eps , x) - \alpha_\eps}{ d(x, \Lambda_\eps)}~,
\end{equation}
where $d(x, \Lambda_\eps)$ is the distance between $x$ and $\Lambda_\eps$ and  
$$ [ \alpha_\eps < h(y_\eps, \cdot )] = \{ ~x \in \R^n ~| ~\alpha_\eps < h(y_\eps, x ) \}~.$$ 
This so called Hoffman constant can be defined following for instance Azé and Corvellec \cite{AC}. Therefore
$$ h(y_\eps , \tilde x) - \alpha_\eps \ge d( \tilde x ,\Lambda_\eps) \, \sigma_{\alpha_\eps} (h) ~.$$
As $d(\tilde x, \Lambda_\eps) = d(\tilde x, \tilde x_\eps ) $ we obtain 
$$  d(\tilde x, \tilde x_\eps ) \leq \frac{h(y_\eps , \tilde x) - \alpha_\eps}{ \sigma_{\alpha_\eps} (h) }~.$$
We have to estimate $ \sigma_{\alpha_\eps} (h) $. In particular we look for $\sigma_o >0$ such that 
$$ \forall \eps \qquad  \sigma_{\alpha_\eps} (h)  \ge \sigma_o~.$$
In \cite{AC}, it is shown that 
$$  \sigma_{\alpha_\eps} (h)  \ge \inf_{h(y_\eps, x) = \alpha_\eps } |\nabla_x h(y_\eps, x)|~,$$
where  $|\nabla_x h(y_\eps, x)|$ stands for the strong  slope of $h$ at $(y_\eps, x)$ with respect to $x$ (\cite{AC}); the strong-slope of a function $ \varphi$ at $x$   is defined as
 $$ |\nabla \varphi (x) | := \left \{ \begin{array}{ll} 
 0 &\mbox{if } x \mbox{ is a local minimum of }\varphi~, \\
 \displaystyle{\limsup_{y \to x } \frac{\varphi(x)-\varphi(y)}{d(x,y)}}&\mbox{otherwise}\end{array}\right.$$
Using (\ref{hyp1}), we have
$$  d(\tilde x, \tilde x_\eps ) \leq \frac{h(y_\eps , \tilde x) - \alpha_\eps}{ \sigma_o  }=\frac{h(y_\eps , \tilde x) - h(y_\eps ,   x_\eps) +o(\eps)}{ \sigma_o   }\to 0~.$$
Indeed $y_\eps \to \bar y $, $x_\eps \to \bar x $, $h$ is continuous and $h(\bar y, \bar x) = h(\bar y , \tilde x)$. 
\\
We may now end the proof.  We can use relation (\ref{eq34}) with $ \tilde x_\eps$ instead of $\tilde x$ so that 
$$\forall z \in C  \qquad h(y_\eps, x_\eps  ) + 2 \eps \, f^2( y_\eps, x_\eps) \le  h(y_\eps, z ) + \eps \, f^2( y_\eps,z) + \eps \, f^2( y_\eps,\tilde x_\eps)   + o(\eps)~;$$
we choose $z= x_\eps$ once again to get
$$    f^2( y_\eps, x_\eps) \le    f^2( y_\eps,\tilde x_\eps)   +\frac{ o(\eps)}{\eps}~.$$
Passing to the limit as $\eps \to 0$ gives  (for every $\tilde x \in \cS(\bar y)$
$$  f^2( \bar y , \bar x) \le    f^2(\bar y, \tilde x )   ~.$$
This means that $\bar x   \in \wcS(\bar y)$. \findemo  

\begin{rem}
It is clear that assumption $(\ref{hyp1})$ is satisfied if $h$ is linear (``linear'' case). Next problem is to find simple conditions for  $(\bar y, \bar x )$ to get  $(\ref{hyp1})$ when $h$ is not linear. One hint is to assume that
   $h$ is $\mathcal{C}^1$ and that $\| \nabla_x h (\bar y, \bar x \|  \neq 0$;  then the strong slope $|\nabla_x h(y_\eps, x)|$  coincides with  the norm  $\| \nabla_xh  (y_\eps, x)\|   $ of the gradient of $h$ with respect to $x$. 
With the convergence  of $(y_\eps, x_\eps)$ to $(\bar y, \bar x )$ (up to a subsequence), there exist $\eps_o$ and $\eta >0$  such that 
$$ \forall \eps \le \eps_o\qquad \|\nabla_xh (y_\eps, x_\eps)\|  \ge \eta > 0~;$$
 next we have to prove that $\|\nabla_xh (y_\eps, x)\|  \ge \eta $ for any $x$ such that 
 $ h(y_\eps,x) = \alpha_\eps$. A good tool could be an ``local inversion theorem'' for the multivalued case but it is not obvious. The problem is still open. We have the same challenge in next section.
 \end{rem}

  %===================================================
 \subsection{Comparison of $(\cP)$ and $(\wcP)$}
%==================================================
Now, it is clear that a solution of the penalized problem $(\cP_\eps)$  is a good approximation of a solution of $(\wcP)$. Anyway, it is not a solution (a priori) of the problem in consideration $(\cP)$. So we have to compare $(\cP)$ and $(\wcP)$. 

The second level of $(\wcP)$ clearly disappears when   the initial problem  lower level    solutions set corresponds to the same revenue for each value of $y$ (or are unique). In this case $(\cP)$ and $(\wcP)$ are equivalent. In other cases, the solution of $(\wcP)$ corresponds to some ``optimal worst'' case solution.

This solution is still important for the decision makers of the upper level problem.
\begin{rem}  Using the same regularization technique, if we replace $\eps f^2(y,z)$ by $-\eps f^2(y,z)$ in the definition of
$h_\eps$, we will obtain (at the limit) an optimal solution of $(\cP)$ which corresponds to an optimal best case solution of our asymmetric game . \end {rem}

    %===================================================
 \Section{Error estimates}
    %===================================================
The purpose of  this section is to study the behavior of $ f(y^*,x^*) -  f(y_\eps, x_\eps) $ as $\eps \to 0$ and provide (if possible) some error estimates.   
Since the penalized problems are nonconvex, we can not use any classical perturbation result.  
We proceed in two steps : we first prove some monotonicity results for the upper level objective function values and then consider classical perturbation analysis results for some auxilliary convex problems. 

\subsection{Preliminary results}

\begin{lemme}\label{lem41} Let  be $\eps>\eps'>0$ and  $y\in K$.   Let be $x_\eps \in \cS_\eps(y )$ and  $\tilde x \in \cS_{\eps'}(y) $.
Then we get
$$  f^2 (y ,x_\eps ) \le f^2 (y , \tilde x ) ~.$$ 
\end{lemme}
\preuve Let us fix  $\eps>\eps'>0$ and choose some $y \in K$.   Let be $x_\eps \in \cS_\eps(y )$ and $ \tilde x \in \cS_{\eps'}(y) $.
Assume that 
 \begin{equation}\label{eq1lam41}
 f^2(y, \tilde x )  <   f^2(y,x_\eps ) ~.
\end{equation}
 
As $ \tilde x \in \cS_{\eps'}(y) $ and $x_\eps \in C$, we have
$$ 
  h(y, \tilde x ) + \eps' f^2(y,\tilde x  )  \le   h(y, x_\eps ) + \eps' f^2(y,x_\eps ) ~,
$$
$$
   h(y, \tilde x ) + \eps' f^2(y,\tilde x  ) + (\eps -\eps')  f^2(y,\tilde x  )   \le   h(y, x_\eps ) + \eps' f^2(y,x_\eps )    + (\eps -\eps')  f^2(y,\tilde x  ) $$
With  (\ref{eq1lam41})  and $\eps>\eps'>0$,  we obtain 
 $$  
   h(y, \tilde x ) + \eps f^2(y,\tilde x  )   <   h(y, x_\eps ) + \eps' f^2(y,x_\eps )  + (\eps -\eps')    f^2(y,x_\eps ) < = h(y, x_\eps ) + \eps f^2(y,x_\eps ) $$
 So 
 $$  
   h(y, \tilde x ) + \eps f^2(y,\tilde x  ) < \mbox{min }  \{ ~ h(y,  x ) + \eps f^2(y, x  ),~x \in C \}$$
   and we get a contradiction.
 \findemo
 %================================
 \begin{lemme}\label{lem42} Let  be $\eps>\eps'>0$ and  $y_\eps $ (respectively $y_{\eps'}$) a solution to $(\cP_\eps)$ (respectively $(\cP_{\eps'})$).   Let be $x_\eps \in \cS_\eps(y_\eps )$ and  $ x_{\eps'} \in \cS_{\eps'}(y_{\eps'}) $.
Then 
$$  f^2(y_\eps ,x_\eps ) \le f^2(y_{\eps'} ,  x_{\eps'}) \leq f^2(y^*, x^*) ~,$$ 
where $y^*$ is a solution to $(\widetilde{\mathcal{P}})$ with $x^*\in \mathcal{S}(y^*)$. 
\end{lemme}
\preuve Using  Lemma \ref{lem41} with  $y=y_\eps$    and $x_\eps \in \cS_\eps(y_\eps )$ gives 
\begin{equation} \forall \tilde x \in \cS_{\eps'}(y_{\eps} ) \qquad 
 f^2(y_\eps ,x_\eps ) \le f^2(y_\eps , \tilde x ) ~.\end{equation}
As  $y_{\eps'}$ is  a solution of $(\cP_{\eps'})$ we get
$$\forall y \in K,~ \forall   x \in \cS_{\eps'}(y)\qquad  f(y_{\eps'}, x_{\eps'})  \ge f(y, x) ~.$$
We may choose in particular $y= y_\eps $ and $x = \tilde x \in  \cS_{\eps'}(y_{\eps})$ to get 
\begin{equation} \forall \tilde x \in    \cS_{\eps'}(y_{\eps}) \qquad f(y_{\eps'}, x_{\eps'})  \ge f(y_\eps,\tilde x)~.\end{equation} 
As $f$ is assumed to be nonnegative we finally obtain
$$  f(y_\eps ,x_\eps ) \le f(y_\eps , \tilde x ) \le  f(y_{\eps'}, x_{\eps'})  ~.$$
Therefore the family $(f(y_\eps ,x_\eps )$ is increasing). The convergence of $f(y_\eps ,x_\eps )$ to $f(y^*,x^*)$  ($f$ is continuous) achieves the proof since $f(y^*,x^*)$ is the limit and the upper bound of the family $(f(y_\eps ,x_\eps ))$ .
\findemo
%============================================

\begin{lemme}\label{lemma43}
Let be $\eps >0$ and $ \tilde x_\eps \in \cS_\eps(y^*)$ where $y^*$ is a solution to $(\widetilde{\mathcal{P}})$.  Then
\begin{equation}\label{eqlem43}
\forall x_\eps \in \cS_\eps(y_\eps) \qquad f(y^*, \tilde x_\eps ) \le f(y_\eps, x_\eps) \leq f(y^*,x^*)~.
\end{equation}  
\end{lemme}
\preuve This is a direct consequence of Lemma \ref{lem42} : the relation $ f(y_\eps ,x_\eps ) \le f(y^*,x^*)$ is obvious
and the relation $ f(y^*, \tilde x_\eps ) \le f(y_\eps, x_\eps) $  comes from the fact that   $y_\eps $ is a solution to $(\cP_\eps)$.
\findemo
\begin{rem} 
 The previous lemmas show that it is sufficient to study 
 $ f(y^*,x^*) -  f(y^*, \tilde x_\eps) $ for some $ \tilde x_\eps \in \cS_\eps(y^*)$.
 \end{rem}
For a large class of realistic problems, the lower level is linearly constrained. moreover, we use 
for our analysis some local error bounds and these bounds are very complicated in case of nonlinear constraints. So, we assume from now on  that  $C$ is polyhedral :
 $$ C = \{ ~x \in \R^n ~\mid ~A x = b ,~x \ge 0~\}~,$$
 where $ A$ is a $m \times n $ real matrix and $b \in \R^m$.  
 
 In the sequel $y^*$ is a solution to   $(\widetilde{\mathcal{P}})$ (which existence is given by Theorem 3.1) and $x^* \in \wcS(y^*) $ ( see (\ref{stildey})) so that 
 $$ x^* \in \mbox{ argmin }\{~f^2(y^*,z)~| ~z \in \mbox{ argmin }\{ h(y^*, \zeta)~, \zeta \in C \}~\} ~.$$
  Let us denote 
  \begin{equation}
\alpha^* = h(y^*, x^*) \mbox{  and } \beta^* = f(y^*,x^*)~.
\end{equation}
Note that $\beta^*$ is the optimal value for $(\widetilde{\mathcal{P}})$ (the upper level) so that we may assume that $\beta^* \neq 0$ (otherwise the problem is trivial). 
We set 
\begin{equation}\label{cstar}
C^* = \{~x \in C ~\mid ~h(y^*,x) \le \alpha^* \mbox{ and } f(y^*, x) \le \beta^*~\}~.
\end{equation}
 Let us give an important property of $C^*$ :
\begin{prop} Assume  $y^*$ is a solution to   $(\widetilde{\mathcal{P}})$  and $C^*$ is defined with (\ref{cstar}), then 
$$C^* = \{~x \in C ~\mid ~h(y^*,x) =\alpha^* \mbox{ and } f(y^*, x) = \beta^*~\}~$$
and
$$C^* = \{~x \in C ~\mid ~h(y^*,x) + f(y^*, x) \leq \sigma^* \stackrel{def}{:=} \alpha^* +\beta^*~\}~.$$ \end{prop}
\preuve Note that it impossible to get $h(y^*,x) \le \alpha^* $, if $x \in C^*$. Indeed, as  $x^* \in \wcS(y^*) $ then $x^* \in S(y^*) = \mbox{ argmin }\{ h(y^*, \zeta)~, \zeta \in C \}~.$ Therefore :
\begin{equation} \label{partiel1}
\forall \zeta \in C \qquad  h(y^*, x^*)  \leq  h(y^*, \zeta)~.\end{equation}
Setting $\zeta=x \in C^*$ gives
$$ \alpha^* = h(y^*,x^*) \leq h(y^*,x)  \leq  \alpha^* ~.$$
So 
$$\forall x \in C^* \qquad h(y^*,x) = \alpha^*~.$$
The same remark holds for $\beta^*$ so that 
\begin{equation}\label{cstar2}                                         
C^* = \{~x \in C ~\mid ~h(y^*,x) = \alpha^* \mbox{ and } f(y^*, x) = \beta^*~\}~.
\end{equation}
Let us call 
$$C' = \{~x \in C ~\mid ~h(y^*,x) + f(y^*, x) \leq \sigma^*  ~\}~.$$ 
It is obvious that $C^* \subset C'$. Conversely, let be  $x\in C'$.  Relation (\ref{partiel1}) yields
$\alpha^* \le h(y^*,x) $ so that 
$$ \alpha^* + f(y^*,x) \leq  \alpha^* + \beta^*~.$$
This gives $ f(y^*,x) \leq   \beta^*$. Similarly, we get $ h(y^*,x) \leq   \alpha^*$ and $x \in C^*$.
\findemo
\\
The main point is now to estimate  (roughly speaking) the distance between the solution $x^*$ and  $\cS_\eps(y^*) $. As $x^* \in C^*$ and $C^*$ is defined with inequalities, we first assume a Hoffman-type condition. 
 %===============================================
  \subsection{ Error estimates under an Hoffman hypothesis}
 %====================================

Following Az\'e and Corvellec \cite{AC} we know that 
$$ \inf_{[\sigma^*< f(y^*,\cdot)+ h(y^*,\cdot) ]} |\nabla_x\left ( f(y^*,\cdot)+ h(y^*,\cdot)\right)|\leq
 \inf_{x\in [\sigma^*< f(y^*,\cdot)+ h(y^*,\cdot) ]}\frac{f(y^*,x)+ h(y^*,x) - \sigma^*}{d(x, [f(y^*,\cdot)+ h(y^*,\cdot)\leq \sigma^*]} ~.$$
 The notation $ [\sigma^*< f(y^*,\cdot)+ h(y^*,\cdot) ]$ stands for the set 
 $$ \{ x \in \R^n~|~\sigma^*< f(y^*,x)+ h(y^*,x) ~\}~.$$
 We note that $[f(y^*,\cdot)+ h(y^*,\cdot)\leq \sigma^*] = C^*$. 
In this subsection, we assume the following 
$$ (\mathcal{H}_1) \qquad   \gamma^* := \inf_{[\sigma^*< f(y^*,\cdot)+ h(y^*,\cdot) ]} |\nabla_x \left( f(y^*,\cdot)+ h(y^*,\cdot)\right )| >0. $$
Let us call $\gamma = \displaystyle{\frac{1}{\gamma^*}}$ : assumption  $(\mathcal{H}_1)$ implies that  
\begin{equation}\label{hoffman1} \forall \varepsilon > 0,~\forall \tilde x_\varepsilon \in \cS_\eps(y^*) \qquad \exists x^*_\varepsilon \in C^*  \mbox{ s.t. } 
\|  \tilde x_\eps - x^*_\eps \| \leq \gamma \left [f(y^*, \tilde x_\eps)+ h(y^*,  \tilde x_\eps) - \alpha^*-\beta^* \right]~.\end{equation}
Note also that 
relation (\ref{eqlem43}) of Lemma \ref{lemma43} yields that 
$$ \forall   \tilde x_\eps \in \cS_\eps(y^*)\qquad f(y^*,  \tilde  x_\eps ) \le \beta^* $$
and $$
h(y^*,  \tilde  x_\eps) \leq \alpha^* + \eps \beta^*$$
because of the definition of $\cS_\eps(y^*)$.Therefore $$ \forall  \tilde  x_\eps \in \cS_\eps(y^*)
\qquad f(y^*,\tilde x_\eps)+ h(y^*,\tilde x_\eps) - \alpha^*-\beta^* \le \eps $$ and
\begin{equation}\label{hoffman2} \
\forall \varepsilon > 0,~\forall \tilde x_\varepsilon \in \cS_\eps(y^*) \qquad \exists x^*_\varepsilon \in C^*  \mbox{ s.t. } 
\|\tilde  x_\eps - x^*_\eps \| \leq  \gamma \eps  ~.\end{equation}
The existence of such Lipschitzian error bound for convex or general inequalities is, itself, an interesting domain of research. It is strongly related to metric regularity properties. A large number of conditions and characterizations can be found in \cite{AC,Az2,Lewis,li,Ng,ye,zhao}. This list of references constitutes a small but significant part of the existent literature.

\begin{rem}  1. Assumption $(\mathcal{H}_1)$ is fulfilled
  if the functions $f$ and $h$ are linear with respect to $x$. Indeed they cannot be identically equal to 0 and the strong slope coincides with the norm of gradient which is a positive constant. \\
  2. $x^*_\varepsilon$ is the projection of  $\tilde x_\eps$ on $C^*$. 
  \end{rem}

\begin{lemme}\label{lem45} Both  $ \tilde x_\eps   \in \cS_\eps(y^*)$ and  $x^*_\eps $ given by (\ref{hoffman2}) converge to $x^* \in   \cS(y^*)$ as $\eps \to 0$.
\end{lemme}
\preuve
We know that $\tilde x_\eps  \to x^*$ (with the previous results).  
Let us set 
 $\displaystyle{ d_\eps = \frac{  \tilde  x_\eps  -  x^*_\eps} {\eps}~.}$
As $ d_\eps $ 
is bounded  (by $\gamma)$ it clear that $x^*_\eps$ and $\tilde x_\eps$ have the same limit point (namely $x^*$).
\findemo

 {In what follows  $\tilde x_\eps$ is an element of $\cS_\eps(y^*)$ and $x^*_\eps $ is the associated element  given by (\ref{hoffman2}) . \\
Let us define 
$$ I(x^*) = \left \{ i \in \{1,\cdots, n\} ~\mid ~x^*_i = 0 ~\right \}~,  \mbox{ and } 
\tilde C = \{ ~d \in \R^n ~\mid ~A d = 0 ~, ~d_{\mid I(x^*)} \ge  0 ~\}~.$$
\textbf{Let $d$ be in  $\tilde C $. }
\\
Then, there exists $\eps_d > 0 $ such that 
 $ \forall \eps < \eps_o , ~ x^*_\eps + \eps d  \in C ~.$ 
 Indeed 
 \begin{itemize}
 \item  $A (x^*_\eps + \eps d )   = A (x^*_\eps)+ \eps Ad  =A (x^*_\eps )    = b~.$
 \item If $i \in  I(x^*) $, then $(x^*_\eps + \eps d)_i \ge x^*_{\eps,i} \geq 0 $.
 \item  If $i \notin  I(x^*) $, then $ x^* _i  >0 $. As $x^*_\eps  \to x^*, \exists \eps_i $ such that 
 $ x^*_{\eps,i} >0$ forall $\eps \leq \eps_i$. Then we choose $\displaystyle{\eta = \inf_{i \notin  I(x^*)  }\{\eps_i\} }$ so that 
 $$ \forall \eps \le \eta \qquad  x^*_{\eps,i} >0~.$$
Now choosing $\eps_d \le \eta$ small enough we get $(x^*_\eps + \eps d)_i \geq 0$ for any $\eps \le \eps_o$.
 \end{itemize}
 As $  \tilde   x_\eps \in \cS_\eps(y^*)$  and  $ x^*_\eps + \eps d \in C$  we have
$$ h(y^*, \tilde  x_\eps) + \eps f^2(y^*,  \tilde  x_\eps)  \leq  h(y^*, x^*_\eps + \eps d  ) + \eps f^2(y^*, x^*_\eps + \eps d  )  ~,
$$
$$
 h(y^*,\tilde  x_\eps) -  h(y^*,   x^*_\eps + \eps d ) + 
  \eps \left [f^2(y^*,  \tilde x_\eps) -   f^2(y^*,   x^*_\eps + \eps d ) \right ]\leq 0 ~.
$$
As the functions are $\mathcal{C}^1$,  we have 
$$  h(y^*,  \tilde x_\eps) =  h(y^*, x^*_\eps) +\nabla_x h(y^*, x^*_\eps)\cdot ( \tilde  x_\eps- x^*_\eps)  + ( \tilde  x_\eps- x^*_\eps) o( \tilde x_\eps- x^*_\eps) $$
\begin{equation}\label{hxeps}
  h(y^*, \tilde  x_\eps) =  h(y^*, x^*_\eps) +\eps \nabla_x h(y^*, x^*_\eps)\cdot d_\eps+ \eps  d_\eps\, o (\eps d_\eps)~,\end{equation}
and
\begin{equation}\label{hxxeps}  h(y^*,   x^*_\eps + \eps d ) =   h(y^*, x^*_\eps) +\eps \nabla_x h(y^*, x^*_\eps)\cdot d +  \eps d \,o (\eps d)~,\end{equation}
where $\nabla_x h$ stands for  the derivative of $h$ with respect to $x$. 
As $ x^*_\eps\in C^*$  and $\tilde  x_\eps \in C$ then 
$$ h(y^*, x^*_\eps) = \alpha^* = h(y^*,x^*) \leq h(y^*,\tilde  x_\eps)~.$$
With relation (\ref{hxeps}) this gives 
$$ \nabla_x h(y^*, x^*_\eps)\cdot d_\eps+   d_\eps\, o (\eps d_\eps) = \frac{ h(y^*,  \tilde x_\eps) -  h(y^*, x^*_\eps)} {\eps} \geq 0~.$$
 As $ d_\eps $ 
is bounded  (by $\gamma)$, there exist cluster points; passing to the limit gives   
\begin{equation}\label{sens1}\nabla_x h(y^*, x^*)\cdot \tilde d  = \lim_{\eps \to 0} \nabla_x h(y^*, x^*_\eps)\cdot d_\eps \ge 0~,\end{equation}
 for  any cluster point $\tilde d$ of the family $d_\eps$.

In addition,  we obtain  with   (\ref{hxeps}) and   (\ref{hxxeps}) $$
 \eps \nabla_x h(y^*, x^*_\eps)\cdot  d_\eps+ \eps  d_\eps \,o (\eps  d_\eps) - \eps \nabla_x h(y^*, x^*_\eps)\cdot d  -  \eps d \,o (\eps d) 
  + 
  \eps \left [f^2(y^*,  \tilde  x_\eps) -   f^2(y^*,   x^*_\eps + \eps d ) \right ]\leq 0 ~,$$
that is 
$$  \nabla_x h(y^*, x^*_\eps)\cdot (  d_\eps - d) +    d_\eps \,o (\eps   d_\eps) - d \,o (\eps d) 
  + 
  \left [f^2(y^*,  \tilde x_\eps) -   f^2(y^*,   x^*_\eps + \eps d ) \right ]\leq 0 ~.$$
Passing to the limit (with Lemma \ref{lem45})  we obtain 
  \begin{equation}\label{genx}
 \nabla_x h(y^*, x^*)\cdot ( \tilde d -d )   \leq 0~,
\end{equation}
where $\tilde d$ is a cluster point of the sequence $d_\eps$ and any $d \in \tilde C$.   
  As $d= 0$ belongs to $\tilde C$, we get 
 $$ \nabla_x h(y^*, x^*)\cdot \tilde d  \le 0~.$$
Finally, we obtain with (\ref{sens1}) 
\begin{equation} \nabla_x h(y^*, x^*)\cdot \tilde d  = \lim_{\eps \to 0} \nabla_x h(y^*, x^*_\eps)\cdot\frac{  \tilde x_\eps -  x^*_\eps}{\eps} = 0~.\end{equation}
 This means that 
 $$\nabla_x h(y^*, x^*_\eps)\cdot( \tilde  x_\eps -  x^*_\eps) = o(\eps).$$ 
 As 
 $$ h ( y^*,  \tilde  x_\eps ) = h(y^*, x^*_\eps) -  \nabla_x h(y^*, x^*_\eps)\cdot  (x^*_\eps -   \tilde x_\eps ) + (x^*_\eps -   \tilde x_\eps )\, o(x^*_\eps -   \tilde x_\eps )$$ we get 
 $$  h ( y^*,  \tilde  x_\eps ) - h(y^*, x^*_\eps)  =  o(\eps) -\eps   d_\eps \, o(\eps  d_\eps ) =  o(\eps) ~.$$
 As $x^*_\eps \in C^*$ then $ h(y^*, x^*_\eps)  = \alpha^*$ and 
\begin{equation}\label{hfin}
\forall \tilde x_\eps \in \cS_\eps(y^*)  \qquad h(y^*, \tilde x_\eps) = h(y^*, x^*) + o(\eps)~.
\end{equation}
As $h$ and $f^2$ play similar roles we have the same result for $f^2$. More precisely

\begin{equation}\label{ffin}
\forall \tilde x_\eps \in \cS_\eps(y^*)  \qquad f^2(y^*,\tilde  x_\eps) - f^2(y^*, x^*) =  o(\eps)~.
\end{equation}
We just proved the following result 
\begin{theo} Assume that $(\mathcal{H}_1) $ is satisfied ; let $y_\eps$ be a solution to  $(\cP_\eps)$  and $\tilde x_\eps \in \cS_\eps (y^*) $.  Then 
$$ h(y^*, \tilde x_\eps)- h(y^*, x^*) = o(\eps)  \mbox{ and }  f^2(y^*,\tilde  x_\eps) - f^2(y^*, x^*) =  o(\eps)\mbox{ as }  \eps \to 0~.$$
Moreover
$$\forall x_\eps \in \cS_\eps (y_\eps) \qquad f(y^*, x^*) - f(y_\eps, x_\eps) = o(\eps) \mbox{ as }  \eps \to 0~.$$

\end{theo} 
\preuve The first assertion has been proved : relations (\ref{hfin}) and  (\ref{ffin}) . We use relation  (\ref{eqlem43}) and the previous result to claim that 
$$ f^2(y^*, x^*)  - f^2(y_\eps, x_\eps) = o(\eps)~.$$
As $f^2(y^*, x^*)  - f^2(y_\eps, x_\eps)  =[ f(y^*, x^*)  + f(y_\eps, x_\eps)]\,[ f(y^*, x^*)  - f(y_\eps, x_\eps)]
$ and $ f(y^*, x^*)  + f(y_\eps, x_\eps) \to 2   f(y^*, x^*)  = 2 \beta^*$ we get the result since  $\beta^* \neq 0$. \findemo

With a bootstrapping technique we obtain the following corollary ;
\begin{cor} Under the assumptions and notations  of the previous theorem, we get  $\forall x_\eps \in \cS_\eps(y_\eps) $ 
$$ \forall n \in \N \qquad  f(y^*, x^*) - f(y_\eps, x_\eps) = o(\eps^n) $$
 and  $\forall \tilde x_\eps \in \cS_\eps(y^*) $ 
 $$h(y^*, \tilde x_\eps) - h(y^*, x^*) = o(\eps^n)~.$$
\end{cor}
\preuve Using relations (\ref{hfin}) and (\ref{ffin}) in assumption $(\mathcal{H}_1) $ we see that relation  (\ref{hoffman2}) becomes 
\begin{equation} 
\forall \varepsilon > 0,~\forall \tilde x_\varepsilon \in \cS_\eps(y^*) \qquad \exists x^*_\varepsilon \in C^*  \mbox{ s.t. } 
\|\tilde  x_\eps - x^*_\eps \| \leq  \gamma o(\eps)  ~.\end{equation}
Using the same technique leads to relations (\ref{hfin}) and (\ref{ffin}) with $\eps^2$ instead of $\eps$ and so on.
\findemo
\begin{rem}
These error estimates are still valid (with the same proofs ) when the penalized problems are approximatively solved.
\end{rem}

 %===============================================
  \subsection{ Error estimates under a ``second-order'' assumption }
 %====================================
If assumption  $(\mathcal{H}_1) $ is not ensured, one may, however, give estimates using the following hypothesis
$$ (\mathcal{H}_2) \qquad \left \{ \begin{array}{c} 
\exists \varepsilon_o > 0~, \exists \delta > 0 , ~\mbox{ such that } \forall x \in  C^* + \mathcal{B}(0, \varepsilon_o)\\
\exists \tilde x \in C^* \mbox{ such that } \| x  - \tilde x  \| ^2\le\delta 
\left [ \left (h(y^*,x) -\alpha^*\right) ^+  +\left (f(y^*, x) - \beta^* \right)^+\right ]
\end{array} \right. $$
\begin{rem}
$ (\mathcal{H}_2)$ means that $C^*$ is H-metrically regular (of the second order). (See the definition of this regularity property for example in \cite{Aus} Def. 4.3.2). $ (\mathcal{H}_2)$ also corresponds to a quadratic growth condition \cite{Bon} Def.3.1 . \\
This assumption is 
significantly weaker than $ (\mathcal{H}_1)$ and covers a large class of problems since it is satisfied when $h(y^*,.) +f(y^*,.) $ is linear or quadratic.
\end{rem}
We have a rather similar result which proof is the the same as in the previous subsection (so that we do not detail it) :
 
 \begin{theo} Assume that $(\mathcal{H}_2) $ is satisfied ; let $y_\eps$ be a solution to  $(\cP_\eps)$ and  $ x_\eps \in \cS_\eps(y_\eps) $.  Then 
$$ f(y^*, x^*) - f(y_\eps, x_\eps) = o(\sqrt{\eps} ) \mbox{ as }  \eps \to 0~,\mbox{ so that}$$
$$\forall \tau >0 \qquad    f(y^*, x^*) - f(y_\eps, x_\eps) = o(\eps^{1-\tau})~.$$
Moreover,  $\forall \tilde x_\eps \in \cS_\eps(y^*) $ 
$$\forall \tau >0 \qquad   h(y^*, \tilde x_\eps) - h(y^*, x^*) =  o(\eps^{1-\tau}) ~.$$

\end{theo} 
\Section{Conclusion} 
With this new penalization approach, the pessimistic formulation of  general bi-level problems becomes, in some way, tractable. Indeed, instead of  the complicated limit problem, we only need to  solve approximately  the penalized one for a small value of the parameter $\varepsilon$. We have given error estimates that  prove that this approximation is reasonable even when Hoffman's  assumption is not satisfied . \\
 \bibliographystyle{plain}

   \end{document}